\documentclass[14pt]{extarticle}

\usepackage{amsmath, amsfonts, amsthm, amssymb}
\usepackage[margin=1in,marginparsep=22pt]{geometry}
\usepackage{hyperref}
\usepackage{mathtools}
\usepackage{enumitem}
\usepackage{float}
\usepackage{graphicx}
\usepackage[symbol]{footmisc}

\newtheorem*{theorem*}{Theorem}
\newtheorem{theorem}{Theorem}
\newtheorem{lemma}{Lemma} 
\newtheorem{proposition}{Proposition}
\newtheorem{corollary}{Corollary}[theorem]

\newtheorem{problem}[theorem]{Problem}

\newtheorem*{conjecture}{Conjecture}

\theoremstyle{definition}
\newtheorem*{definition}{Definition}

\newtagform{brackets}[\textbf]{$\langle$}{$\rangle$}
\usetagform{brackets}

\newcommand{\tscaption}[1]{\caption{\textsc{\small #1}}}

\let\gte\geqslant
\let\lte\leqslant

\newcommand\vD{\mathcal{D}}

\newcommand\N{\mathbb{N}}
\newcommand\Z{\mathbb{Z}}
\newcommand\Zp{\Z_{\gte 0}}
\newcommand\GP{\mathbb{GP}}

\def\ifempty#1{\def\temp{#1}\ifx\temp\empty}
\newcommand\tauapprox{1.8637}
\newcommand\broot[2]{\root{\raise1pt\hbox{$\scriptstyle#1$}}\of{#2}}
\newcommand\bslash{\kern-1.1pt\setminus\kern-2pt}
\newcommand\ovrl[1]{\overline{\vphantom{\hbox{\rule{0pt}{7.5pt}}}#1}}

\newcommand{\pic}[2]{\includegraphics[scale=#1]{#2}}

\newcommand{\clpic}[4]{\begin{figure}[H]%
\begin{center}\noindent\pic{#1}{#2}\end{center}%
\ifempty{#4}\empty\else\label{fig:#3}\tscaption{#4}\fi%
\end{figure}}
\newcommand{\dlpic}[6][0pt]{\begin{figure}[H]%
\begin{center}\pic{#2}{#3}\hspace{#1}\pic{#2}{#4}\end{center}%
\ifempty{#6}\empty\else\label{fig:#5}\tscaption{#6}\fi%
\end{figure}}

\allowdisplaybreaks

\begin{document}

\title{Elementary proof for the bounds of the complexity of a planar multigraph and the size of a prime rectangular squaring}
\author{Dmitri Fomin}
\date{\today}

\maketitle

\section{Summary}

This article contains two results together with their relatively elementary proofs. 

The first one (see Theorem \ref{thm:PlanarGraphSpanTrees}) presents the upper boundary on the number of spanning trees in a finite planar multigraph.

\begin{theorem*}[Upper bound for the complexity of a planar multigraph]
The complexity (the number of spanning trees) of a planar multigraph with $n$ edges does not exceed $\tau^n$, where $\tau \approx \tauapprox$.
\end{theorem*}

This result is, quite possibly, already known and/or published---my quick web search did not turn up anything but that does not really prove much. It also seems plausible that this inequality is actually true for the ``best possible'' value of $\tau^* \approx 1.7916$ (see Conjecture \ref{thm:CatalanNumberConjecture}).

The second result (see Corollary \ref{thm:SizeOfPrimeSquaring}) uses the above theorem to improve on the well-known Conway's inequality for the number of tiles in a prime rectangular squaring (see \cite{Conway}).
 
\begin{theorem*}[Lower bound for the size of a rectangular squaring]
If $M \times N$ rectangle ($M \lte N$) is dissected into $n$ squares in such a way that the greatest common divisor of their sizes is $1$, then $n \gte \log_\tau N \approx 1.1134 \log_2 N$.
\end{theorem*}

\section{Motivation}

The original motivation for this article came from the following relatively simple contest-style question (see item (b) of the following problem).

\begin{problem}
a) Square $S$ with dimensions $N \times N$ consists of $N^2$ unit squares. A corner unit square is removed. Prove that if $N = 2^n$, then the remaining part cannot be dissected into less than $n$ squares.

b) An arbitrary unit square (not necessarily the corner one) is removed from $S$. Prove that the remaining part still cannot be dissected into less than $n$ squares.

c)* Prove some better lower bound such as, perhaps, $4n/3$. (It is quite easy to see that in both (a) and (b) the remaining part can always be dissected into $3n$ squares.)
\end{problem}

It can be easily shown (at high school math level) that in part \textit{(a)} if number $k$ is such that $N \gte \Phi_k$, where $(\Phi_k)$ are the Fibonacci numbers $\Phi_1 = 1$, $\Phi_2 = 1$, $\Phi_n = \Phi_{n-1} + \Phi_{n-2}$, then the squaring must consist of at least $k+1$ squares, which gives us a noticeably better result than the one stated in part \textit{(a)}, since $\Phi_k$ is the integer closest to $\phi^k/\sqrt{5}$, where $\phi = (1+\sqrt5)/2$ is the golden ratio. However, that proof only works when the said unit tile lies in the corner.

Our investigation of parts \textit{(b)} and \textit{(c)} proceeds in the obvious direction of the well developed theory of rectangular squarings, which employs graph theory as well as the theory of electric circuits. An additional motivation is to present all proofs in a purely graph-theoretic way without any use of somewhat vague ``electricity''-based reasoning, which often does look quite informal and not entirely convincing to a ``discrete'' mathematician (we will use some electricity-related facts and terminology here but only to point out some general analogies as well as to underline the motivation for the terms and ideas).

\section{Lemmas}

The first supporting lemma is the following useful fact, which I was not yet able to find in the published literature---however it seems rather unlikely that this (or some similar and possibly better result) was not known before.

\begin{theorem}
\label{thm:PlanarGraphSpanTrees}
For any finite planar multigraph $G$ with $n$ edges, the number of its spanning trees $t(G)$ does not exceed $\tau^n$, where $\tau = 1.8637065...$ is the largest root of cubic equation $x^3 - x^2 - 3 = 0$.
\end{theorem}

\begin{proof}
Let us define function $\mu:\Zp \rightarrow \Zp$ as
\begin{equation}
\label{eqn:defMu}
\mu(n) = \max \{t(G): G \in \GP_n \} \,,
\end{equation}
where $\GP_n$ is the set of all planar multigraphs with $n$ edges.

We will prove inequality $\mu(n) \lte \tau^n$ by induction on $n$. The basis of induction immediately follows from observing that $\mu(0) = \mu(1) = 1$, $\mu(2) = 2 < \tau^2$, $\mu(3) = 3 < \tau^3$.

Obviously, we can assume $G$ to be connected. Indeed, if $G = \bigcup G_k$, where $G_k$ are multigraph $G$'s components of connectedness, then $t(G) = \prod t(G_k)$.

Any multigraph $G \in \GP_n$ satisfies the well-known recursive formula
\begin{equation}
\label{eqn:tGe}
t(G) = t(G\bslash e) + t(G/e) \,,
\end{equation}
where $e$ is an arbitrary edge in $G$, $G\bslash e$ is the multigraph obtained from $G$ by removing edge $e$, and $G/e$ is the multigraph obtained by contracting edge $e$ (the proof is quite simple: the set of all spanning trees in $G$ can be split into two disjoint subsets---one containing all spanning trees with edge $e$, and the other one comprising all spanning trees without edge $e$).

Note that both graphs $G\bslash e$ and $G/e$ belong to $\GP_{n-1}$.

\textbf{Case 1.} $G$ contains cycle $C$ of length not exceeding $3$. 

If cycle $C$ has length $1$ (it is a loop) then 
$$
t(G) = t(G\bslash C) < \mu(n-1) < \tau^{n-1} < \tau^n \,.
$$

If $C$ has length $2$ (double edge), then $C$ consists of two edges $e'$ and $e''$ connecting the same two vertices $u, v \in G$. Then at least one of these two edges must be excluded from any spanning tree and therefore we have
$$
t(G) = t(G\bslash \{e', e''\}) + 2t(G/e'/e'') < 3\mu(n-2) < 3\tau^{n-2} < \tau^n \,,
$$
since $\tau > \sqrt{3} \approx 1.7321$.

Last subcase: cycle $C$ has length $3$, $C = (uvw)$.

Then by using \eqref{eqn:tGe} we can obtain
\begin{align*}
&t(G) = t(G\bslash (uv)) + t(G/(uv)) \,, \\
&t(G/(uv)) = 2 t(G/(uv)\bslash \{(uw), (vw)\}) + t(G/(uv)/(uw)/(vw)) 
\end{align*}
and therefore
\begin{align*}
t(G) &= t(G\bslash (uv)) + 2 t(G/(uv)\bslash \{(uw), (vw)\}) + t(G/(uv)/(uw)/(vw)) \\
&\lte \mu(n-1) + 3 \mu(n-3) \lte \tau^{n-1} + 3\tau^{n-3} = \tau^n \,,
\end{align*}
since $\tau^3 - \tau^2 - 3 = 0$.

\clpic{1.0}{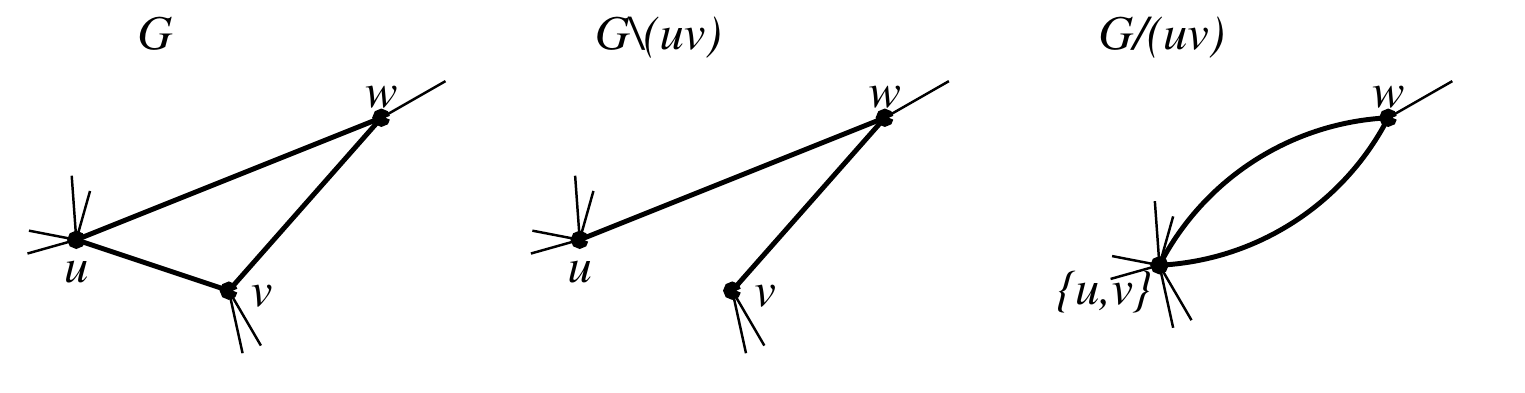}{EdgeOperations}{Deletion and contraction of an edge (Case 1, cycle length 3)}

Now consider planar graph $G^*$ dual to $G$. Then $t(G) = t(G^*)$ (the proof is quite easy as the duality provides a one-to-one correspondence between spanning trees in $G$ and $G^*$; simply put, for any spanning tree $T \in G$ the subgraph in $G^*$ formed by all edges that do not intersect edges of $T$ is a spanning tree for $G^*$). Also note that the number of edges in $G$ is the same as the number of edges in $G^*$. Finally, the degree of any vertex in $G^*$ equals the number of boundary edges of the corresponding planar face of $G$, and vice versa.

It follows that if any vertex in $G$ has degree $3$ or less, we can apply \textbf{Case 1} to dual graph~$G^*$ to prove the induction step for both $G^*$ and~$G$.

\textbf{Case 2.} All vertices in $G$ have degrees greater than or equal to $4$ and all planar faces of $G$ (including the exterior face) have at least $4$ boundary edges.

In this case, let $m$ be the number of vertices in $G$, and $k$---the number of planar faces in $G$. Then obviously we have $2n \gte 4m$ and $2n \gte 4k$, which contradicts the Euler formula for planar graphs $n + 2 = m + k$. Therefore such a planar multigraph does not exist.

This completes our proof of Theorem \ref{thm:PlanarGraphSpanTrees}.
\end{proof}

Most likely, Theorem \ref{thm:PlanarGraphSpanTrees} is valid for $\tau^* = e^{2C/\pi} \approx 1.7916$, where 
$$
C = 1 - \frac{1}{3^2} +\frac{1}{5^2} - \frac{1}{7^2} + \ldots 
= 0.9159655941772\ldots \,,
$$
the so-called Catalan constant. If that were indeed so, then $\tau^*$ would be the tight upper bound, since it is well known that for graph $R_k$ (the rectangular $k \times k$ planar grid), which has $n = k(k+1)$ edges, we have
$$
\lim_{n \to \infty} \frac 1n \ln(t(R_k)) = \frac {2C}{\pi} \,,
$$
see \cite{Wu}.

\begin{conjecture}
\label{thm:CatalanNumberConjecture}
$$
\mu(n) \lte {\mathcal C}^n \,, \text{ where } \mathcal C = e^{2C/\pi} \,.
$$
(function $\mu(n)$ is defined in \eqref{eqn:defMu} as the maximum number of spanning trees for a planar multigraph with $n$ edges).
\end{conjecture}

\smallskip

The second lemma required for our investigation can be called the ``edgewise'' analog of the famous Kirchhoff's matrix tree theorem (actually, the matrix we consider in that lemma is in fact the matrix that Kirchhoff himself originally investigated, see \cite{Kirchhoff}, \cite{Bryant}, and \cite{Kirby}).

\begin{definition}
For any finite connected edge-weighted directed graph $G$ (such a graph is often called a \textit{network}), and for any set $P$ of its edges with assigned orientation (which does not have to be the same they carry in $G$) let us define the \textit{oriented weight} of $P$ as
$$
\ovrl{\omega}(P) = \sum_{(\overrightarrow{uv}) \in P} s_{uv} \cdot \omega(uv),
$$
where $\omega(uv)$ is the weight of edge $uv$ and $s_{uv}$ equals $+1$ or $-1$ depending on the coincidence of the orientation assigned to $(uv)$ with the orientation that $uv$ carries in~$G$.
\end{definition}

\begin{definition}[Edgewise Kirchhoff's matrix of multigraph]
\label{def:EdgewiseKirchhoffMatrix}
For any finite connected indexed multigraph $G$ with $n$ edges and $m$ vertices let us construct the following system of $n$ linear equations in $n$ variables\footnote{We call multigraph with $m$ vertices and $n$ edges \textit{indexed} if its vertices are labeled by indexes $1$ through $m$ and its edges are labeled by indexes $1$ through $n$.}.

Assume that $n > 0$, and assign variable $\omega_{uv}$ to each edge $uv \in G$, treating it as the weight of that edge, thus turning $G$ into a network.

Choose two arbitrary vertices $S, T \in G$ (we will call them the \textit{poles}), and then also choose an arbitrary orientation for all edges in $G$.

The first $m-1$ equations of the system have the form
$$
[1] \qquad \ovrl{\omega}(\{uv\}_{(uv) \in G}) = f_u, 
$$
where vertex $u$ is fixed and different from $S$, the summation is done over the set of all vertices $v$ adjacent to and different from $u$, and numbers $f_u$ are zeros except, possibly, for $f_T$. In other words, these equations correspond to the Kirchhoff's \textbf{current} law of electricity (total current flowing into a non-pole vertex $u$ is equal to the total current flowing out of $u$). In this notation edge $(uv)$ always carries orientation $u \rightarrow v$ (its orientation in $G$ may be different).

The remaining $n-m+1$ equations of the system correspond to the Kirchhoff's \textbf{voltage} law (sum of voltages around any closed loop is zero). Namely, we choose any spanning tree $F$ in multigraph $G$ and consider $n-m+1$ edges constituting the complement $G\bslash F$. Each such edge $e$ defines unique simple cycle $C = (u_1, u_2, \ldots, u_k, u_{k+1} = u_1)$ in $F\cup{e} \subset G$, so we write $n-m+1$ equations of the form
$$
[2] \qquad \ovrl{\omega}(\{u_iu_{i+1}\}_{i = 1, \ldots, k}) = 0.
$$

\begin{figure}[H]
\centerline
{%
\begin{tabular}{ccc}
\pic{1.0}{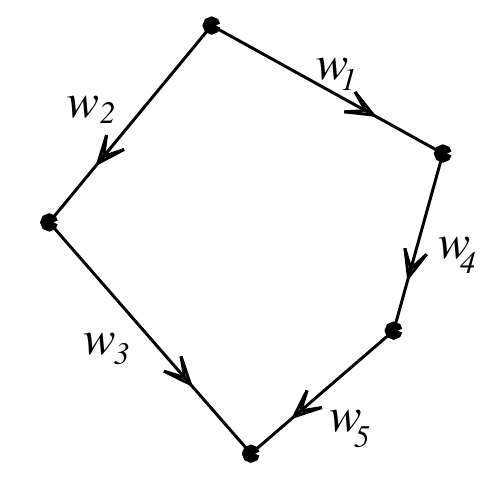} & 
\qquad\qquad\qquad &
\smallskip
\pic{1.0}{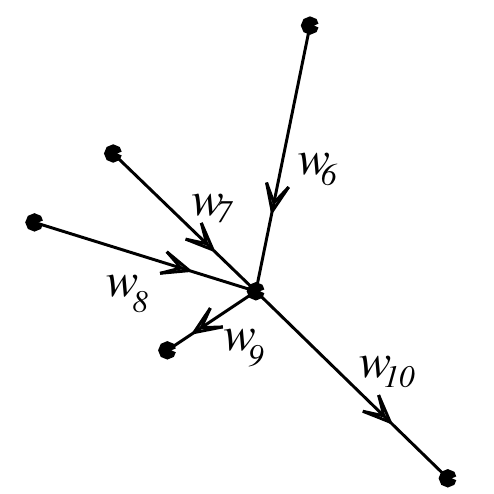} \\
$\omega_1 + \omega_4 + \omega_5 - \omega_2 - \omega_3 = 0$ & &
$\omega_6 + \omega_7 + \omega_8 - \omega_9 - \omega_{10} = 0$
\end{tabular}
}
\label{fig:EdgeEquations}
\tscaption{Kirchhoff's equations for the edge weight variables}
\end{figure}

We will denote the matrix of this system as $M_{G,F}$ (or often simply $M_G$) and call this matrix an \textit{edgewise Kirchhoff's matrix} of multigraph $G$.
\end{definition}

\begin{definition}
Using notation of Definition \ref{def:EdgewiseKirchhoffMatrix}
define the \textit{sign} $\sigma(F)$ for any spanning tree $F$ in indexed network $G$ with some pole vertex $S$ as follows.

Denote indexes of all vertices except $S$ by $a_1 < \ldots < a_{m-1}$. Similarly, edges of $F$ have indexes $b_1 < \ldots < b_{m-1}$.

First step: write down sequence $\alpha$ consisting of $m-1$ zeros (or of the same number of empty placeholders).

Second step: for every vertex $v \in F' = F\bslash \{S\}$ there is exactly one simple path in $F$ connecting $v$ and $S$, and therefore we can define edge $e_v \in F$ as the first edge on that path (it is, of course, incident to $v$).

Edge $e_v$ is either oriented away from or towards $v$ (or equivalently, towards or away from $S$). Label $e_v$ with $1$ or $0$ accordingly.

Find index $a_i$ of vertex $v$ as well as the index $b_j$ of edge~$e_v$, then set $\alpha_i = j$.

Repeat second step for every vertex of $F'$. It is quite clear that after this process is over, sequence $\alpha$ becomes a permutation of $\N_{m-1}$. Also now all edges of tree $F$ are labeled with zeros and ones. 

\clpic{0.9}{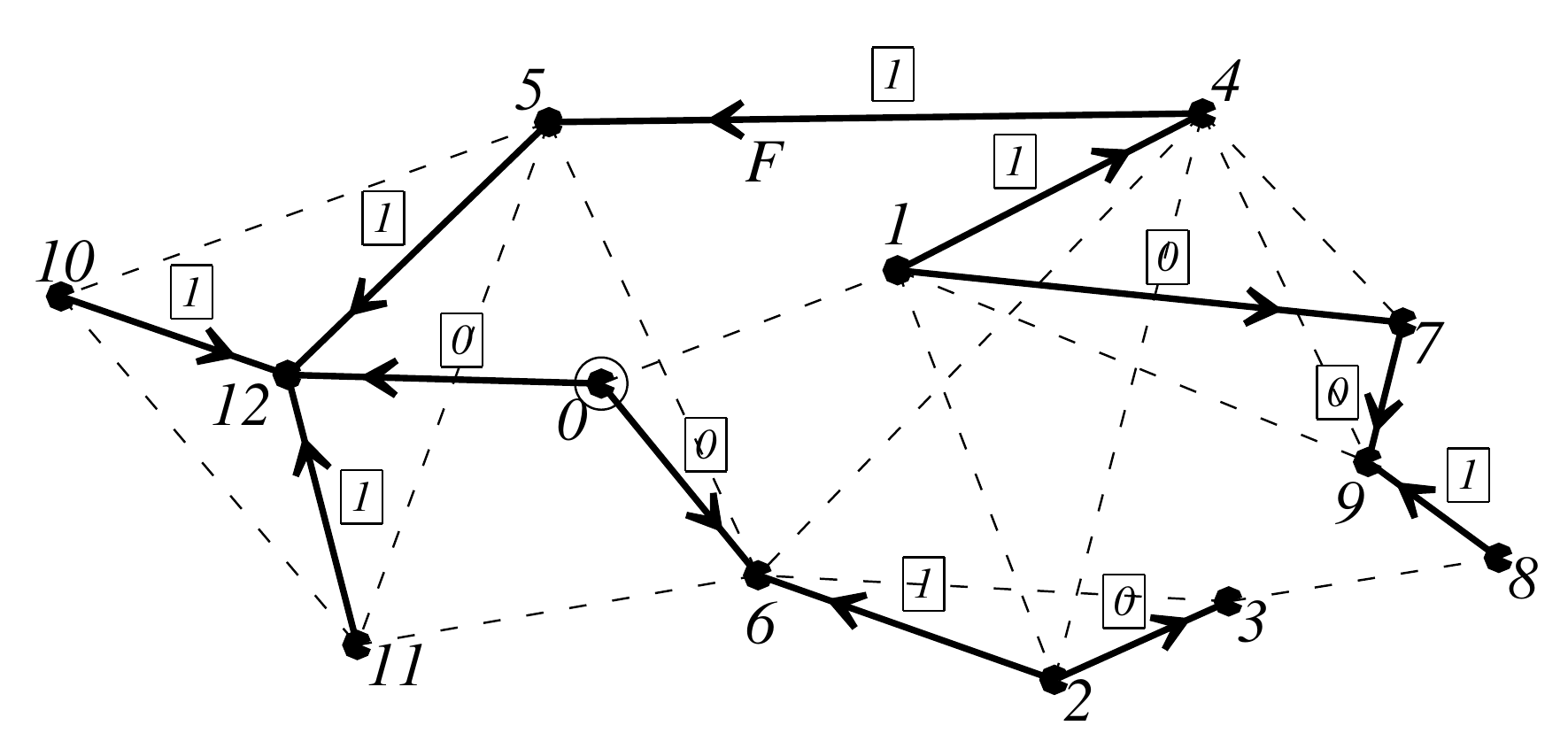}{spanTreeLabeling}{Spanning tree with labels (pole $S$ is circled)}

Now, let $p_1$ be the sum of all the edge labels, $p_2$---the parity of permutation $\alpha$ (that is, $0$ if $\alpha$ is even and $1$ if it is odd), and $p_3$---the parity of the sum 
$\sum_{j = 1}^{m-1}(b_j + j)$.

Finally, define $\sigma(F) = (-1)^{p_1 + p_2 + p_3}$.

\end{definition}

\begin{definition}
Let us denote the Type [1] row of matrix $M$ which corresponds to vertex $u$ by $r_u$, and the Type [2] row which corresponds to edge $e \in G\bslash F$---by $r_e$. Accordingly, for any edge $e \in G$ we will denote by $c_e$ that column of $M$ which corresponds to variable $\omega_e$.

Consider all Type [2] rows $r_e$ and write out---in the same order---the indexes of the corresponding columns $c_e$. Let $\lambda(M)$ be $(-1)^{p_4}$ where $p_4$ is the parity of the resulting pseudo-permutation.\footnote[1]{Any sequence $s$ of $k$ different numbers can be transformed into a permutation of $\N_k$ if we replace each number in it with its index in sequence $\tilde{s}$ which is the result of sorting $s$ in the ascending order.} Of course, if the rows $r_e$ are written out in accordance with the order of edges, ie, in the same order as columns $c_e$, then $p_4 = 0$ and $\lambda(M) = 1$. 

Define $\rho(M) = \sigma(F)\lambda(M) = (-1)^{p_1 + p_2 + p_3 + p_4}$.

Finally we define function $\eta(G)$ by the formula
$$
\eta(G) = \
\begin{cases*}
1                  & \hbox{if $n=0$} \\
\rho(M) \det(M_G)  & \hbox{if $G$ is connected} \\
0                  & \hbox{otherwise}
\end{cases*}
$$
\end{definition}

\begin{proposition}
\label{thm:EtaWellDefined}
Function $\eta(G)$ is well defined, meaning that its value does not depend on the vertex or edge indexing, on the choices of the poles, on the edge orientations, or on spanning tree $F$.
\end{proposition}

\begin{proof}

Let us switch indexes for some two vertices. Then the only component of $\rho(M) = (-1)^{p_1 + p_2 + p_3 + p_4}$ that changes its parity is $p_2$. Since that operation also inverts the sign of $\det(M)$, this proves that the vertex indexing is not relevant.

Now transpose the order of two edges with indexes that differ by one. If they both belong to $F$, the only affected component is $p_2$. If they both are outside of $F$ then the only affected component is $p_4$. If one belongs to $F$ and the other does not, then the only affected component is $p_3$. Therefore in any case there is exactly one component changing its parity, and since the sign of $\det(M)$ also changes, we are done with this part of the proof.

The choice of edge orientations also does not matter. Let us flip the orientation of edge $e$. If $e \in F$ then column $c_e$ of matrix $M$ will be multiplied by $-1$ and so $\det(M_{G,F})$ changes its sign. Accordingly, the only component of $\rho(F)$ that changes its parity is $p_1$. If $e \notin F$ then both column $c_e$ and $r_e$ change signs and $\det(M)$ stays the same while none of $p_i$ are affected.

Furthermore, selecting another spanning tree $F$ does not change the outcome. Indeed, every other spanning tree can be obtained from $F$ by several operations of the following type. Add any edge $e \in G\bslash F$, thus creating cycle $C$, then remove an adjacent edge $e' \in C \cap F$, resulting in spanning tree $F'$. Let us denote original matrix $M_{G,F}$ by $M$, and new matrix $M_{G,F'}$ by $M'$---now we must prove that $\det(M) = \det(M')$.

Switching from $F$ to $F'$ affects only Type [2] (``voltage'') rows in matrix $M$; more precisely, it
\begin{enumerate}[label=(\alph*),itemsep=-1pt]
\item replaces row $r_e$ (defined by cycle $C$) by row $r_{e'}$ which corresponds to the same cycle $C$ in matrix $M'$;
\item changes some coefficients in the rows of $M$ which correspond to cycles that pass through edge $e'$ but not through $e$.
\end{enumerate}

Obviously, step (a) does not change our matrix at all since $r_e$ in $M$ and $r_{e'}$ in $M'$ represent the same equation $\ovrl{\omega}(C) = 0$ in both cases. And the change in any row $r$ from step (b) simply consists of subtracting row $r_e$ from $r$---and that operation does not affect the value of the matrix determinant.

Finally, the choice of the poles is also irrelevant. Obviously, we only need to prove it for vertex~$S$.

Choose some other vertex $S'$ as the pole instead of $S$ and  consider again the corresponding matrices $M$ and $M'$. Clearly, the sum of all Type [1] (``current'') rows in $M$ is equal to
$$
\sigma = \ovrl{\omega}(\{Su\}_{(Su)\in G}) \,.
$$
Therefore replacing row $r_{S'}$: $\ovrl{\omega}(\{S'u\}_{(S'u)\in G})
$ in $M$ with row $r_S$: $\ovrl{\omega}(\{Su\}_{(Su)\in G})$ in $M'$ is equivalent to adding to $r_{S'}$ the sum of all other Type [1] rows of $M$; once again, this operation does not change the value of the determinant. Incidentally, $\sigma = f_T$; we will call this number the \textit{throughput} of network $G$.
\end{proof}

\begin{theorem}[Edgewise Kirchhoff's matrix tree theorem]
\label{thm:EdgewiseKirchhoff}
For any finite connected multigraph $G$ we have $\eta(G) = t(G)$.
\end{theorem}

\begin{proof}

Despite the fact that this theorem was already proved---albeit in a slightly different form---in \cite{Kirchhoff} and in \cite{Bryant}, we will present here another proof, which is a tad more graph-theoretical rather than linear-algebraic.

First, note that the special cases---when $G$ is not connected or when it has zero edges---are trivial.

Second, if $G$ has any loops, they can be eliminated---it is easy to see that this operation affects neither $t(G)$ nor $\eta(G)$.

Let us use the induction on $n$, employing the recursive formulas \eqref{eqn:tGe}. If we prove that $\eta(G)$ satisfies the same recurrent equation, namely that for an arbitrary edge $e \in G$ equality
$$
\eta(G) = \eta(G\bslash e) + \eta(G/e)
$$
holds true, then equality $\eta(G) = t(G)$ will immediately follow by induction on the number of edges (the basis of induction with $n=0, 1$ is quite self-evident).

We will use notation $G' = G\bslash e$ and $G'' = G/e$; we will also assume that edge $e$ connects vertices $u$ and $v$ different from $S$ (the other case is nearly identical).

\textbf{Case 1.} $G'$ is not connected. This means that $\eta(G') = 0$; also it follows that in $M_G$ variable $\omega_e$ is not present in any equation of Type [2] as $e$ cannot belong to any simple cycle in $G$. There are exactly two Type [1] rows---$r_u$ and $r_v$---in which $\omega_e$ is used. If we replace row $r_u$ with $r_u + r_v$, then the determinant of $M$ will not change. In the resulting matrix $M^*$ column $c_e$ has only one non-zero element, namely $\pm 1$ in row $r_v$. Therefore, removing column $c_e$ and row $r_v$ from $M^*$ will produce $(n-1)\times (n-1)$ matrix $M''$ whose determinant equals $\det(M_G) (-1)^d$. It is easy to see that this matrix is $M_{G''}$, $\sigma(F)$ stays the same and $(-1)^d$ coincides with the ratio $\lambda(M_G)/\lambda(M_{G''})$. This proves Case 1.

\textbf{Case 2.} $G'$ is connected. Proposition \ref{thm:EtaWellDefined} allows us to choose any spanning tree $F$ in $G$. Since $G' = G\bslash e$ is connected, we can find such a tree in $G'$ and then set $M_G = M_{G,F}$. We can also assume that column $c_e$ is the last column (index $n$), that row $r_e$ is also the last row of matrix~$M_G$, and that edge $e = \overrightarrow{uv}$.

Column $c_e$ in $M_G$ has three non-zero elements---$1$ in Type [2] row $r_e$, and two more with opposite signs in Type [1] rows $r_u$ and $r_v$ (more precisely, $-1$ in row $r_u$ and $1$ in row $r_v$). Let $k$ be the index of row~$r_v$.

Replace row $r_u$ with $r_u + r_v$ to produce matrix $M^*$ with $\det(M^*) = \det(M_G)$. Column $c_e$ in $M^*$ now has only two non-zero elements---at the intersections with rows $r_e$ and $r_v$. 

Thus we have
$$
\det(M_G) = \det(M^*) = (-1)^{2n}\det(M_1) + (-1)^{n+k}\det(M_2) \,
$$
where $M_1$ and $M_2$ are the minors corresponding to those two non-zero elements in column $c_e$. Multiplying by $\rho(M_G)$ gives us the equality
\begin{equation}
\label{eq:SumOfTwoCofactors}
\eta(G) = \rho(M_G)\det(M_1) + (-1)^{n+k}\rho(M_G)\det(M_2) \,.
\end{equation}

Minor $M_1$ is obtained by removing column $c_e$ and row $r_e$ (recall that $e \in G\setminus F$). If we perform the row operation inverse to the one described above---namely, subtracting row $r_v$ from the former row $r_u$ which is now equal to $r_u + r_v$---then the result will obviously be equal to matrix $M_{G'}$. Therefore $\det(M_1) = \det(M_{G'})$. Also $\rho(M_G) = \rho(M_{G'})$ as parity of neither component $p_i$ is affected by the removal of column $c_e$ and row $r_e$. Thus we conclude that the first summand on the right-hand side of \eqref{eq:SumOfTwoCofactors} equals $\eta(M_{G'})$.

\clpic{1.0}{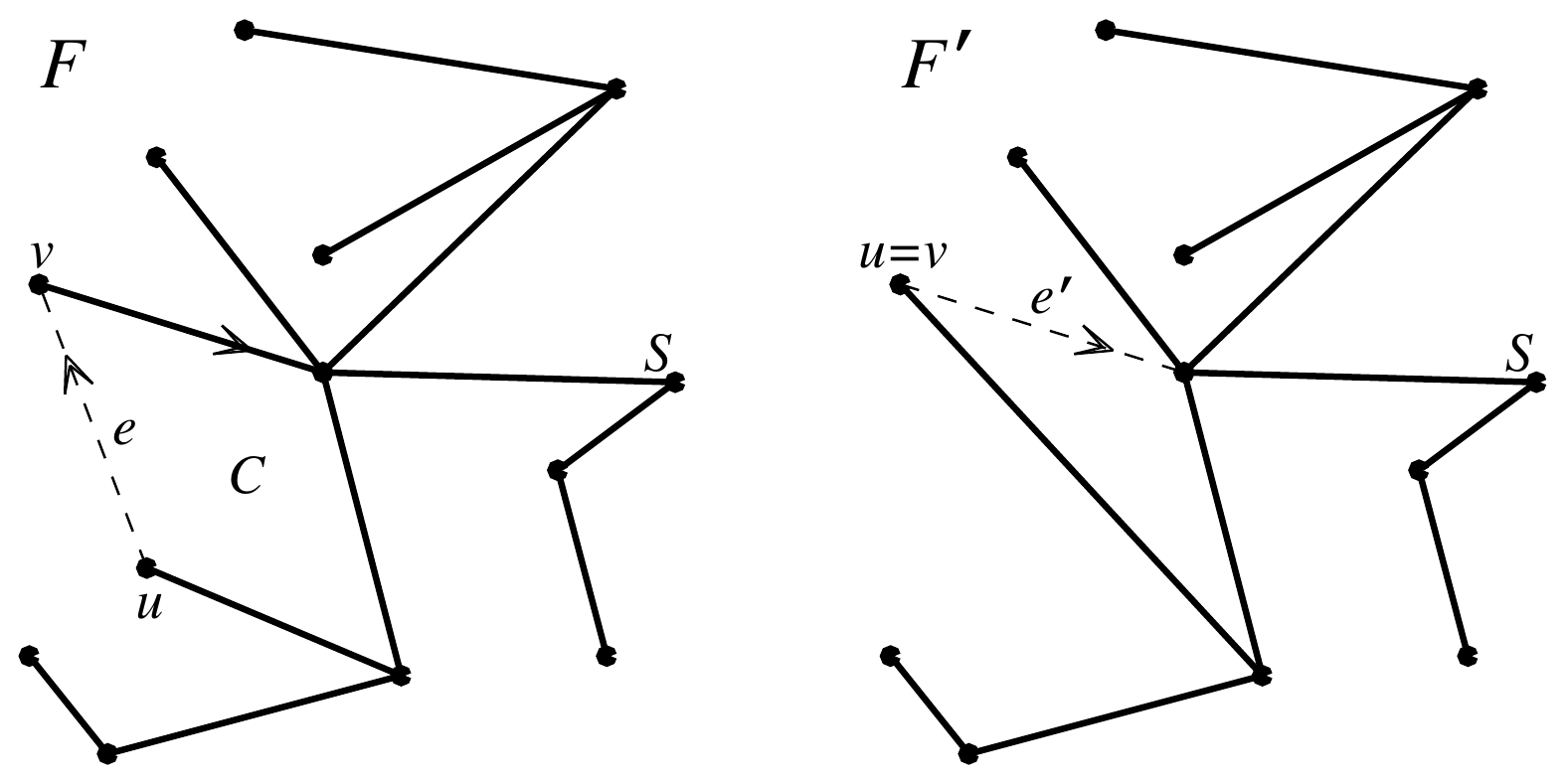}{spanTreeContraction}{Spanning tree ``contraction''}

Similarly, the second minor, which is obtained from $M^*$ by removing column $c_e$ and row $r_v$, is matrix $M_{G''}$. Indeed, since $F/e$ is no longer a tree, in order to produce matrix $M_{G''}$ we use instead $F' = (F/e)\bslash e'$, where $e'$ can be any edge in the simple cycle $C$ in $F/e$ created by contracting edge $e$. Then the only difference between these matrices will consist in replacing row $r_e$ with row $r_{e'}$--but that row is identical to $r_e$ when monomial $\pm w_e$ is removed.

Since we are free to select edge $e'$ in cycle $C$, we choose the one which shared common vertex $v$ with $e$ inside~$C$. Without loss of generality we can assume that $e$ and $e'$ were co-ori\-en\-ted in that cycle. Let $q$ be the index of $e'$, and $d$ is the number of edges in $F'$ with indexes greater than $q$. Then the change in components of $\rho(M_{G''})$ versus components of $\rho(M_G)$ can be quickly expressed as follows. 

Change in $p_1$ is equal to~$1$. Change in $p_2$ equals $d-k \pmod 2$. Change in $p_3$ is $2m-1+(n-q) \equiv n-q-1 \pmod 2$, and finally, change in $p_4$ is $d-q$. Adding these changes together we get $1 + (d-k) + (n-q-1) + (d-q) \equiv n+k \pmod 2$. It follows that
$$
(-1)^{n+k}\rho(M_G) = \rho(M_{G''}) \,.
$$

Therefore the second summand in \eqref{eq:SumOfTwoCofactors} is $\eta(M_{G''})$ and the proof of the recursive formulas for $\eta(G)$ is complete.
\end{proof}

\begin{corollary}
The number of the spanning trees of a finite connected multigraph equals the absolute value of the determinant of its edgewise Kirchhoff's matrix.
\end{corollary}

\section{Main Theorem}

Now we go back to rectangular dissections (or tilings) of integer-sized rectangles, that is, rectangles with dimensions $M\times N$ such that $M, N \in \mathbb N$.

\begin{definition} For any integer rectangle $R$ we will denote its maximum side length by $\max(R)$, its minimum side length---by $\min(R)$, and the greatest common divisor of these two positive integers---by $\gamma(R)$.
\end{definition}

\begin{definition} For any integer rectangle $R$ we will define \textit{arithmetic bias} $\delta(R)$ (or simply \textit{bias}) of rectangle $R$ as $\delta(R) = \max(R)/\gamma(R)$.
\end{definition}

\begin{definition} A rectangular dissection is called a \textit{squaring} if all of its members (tiles) are squares. It is called a \textit{prime} squaring if the greatest common divisor of all tile sizes equals $1$.
\end{definition}

What follows is a relatively straightforward ``generalization'' of Theorem $4$ from John Conway's paper \cite{Conway}, where, among other things, Conway proved that any prime squaring of a rectangle one of whose sides equals $N$ must contain at least $\log_2 N$ tiles. 

\begin{theorem}
\label{thm:SumOfBiases}
For any rectangular dissection $\{R_i\}$ of integer rectangle $R$ consider number $d = \gcd\{\gamma(R_i)\}$ (that is, $d$ is the greatest common divisor of all dimensions of the tiles). Then 
\begin{equation}
\label{eqn:sumofdeltas}
\sum \delta(R_i) \gte \log_\tau \left( \frac{\max(R)}{d} \right) \,,
\end{equation}
where $\tau \approx \tauapprox$ is the number defined in the statement of Theorem \ref{thm:PlanarGraphSpanTrees}.

\end{theorem}

\begin{proof} 

First, let us prove inequality \eqref{eqn:sumofdeltas} for the prime squarings. Since bias of a square tile is $1$, and $\gcd\{\gamma(R_i)\} = 1$, then we can reformulate it as follows.

\begin{lemma}
\label{thm:TauConwayLemma}
For any prime squaring $\vD$ of integer rectangle $R$ the number of squares in dissection cannot be less than $\log_{\tau}(\max(R))$.
\end{lemma}

\begin{proof}

First, let us assume that the largest dimension of $R$ is the horizontal one, and so $\max(R) = N$, where $N$ is the length of the bottom side of $R$ (if the vertical dimension is the larger one, we can simply rotate the rectangle by $90^\circ$).

For the given dissection $\vD$ of $R$ into $n$ square tiles we construct a network\footnote[1]{We remind the reader that a network is a directed connected edge-weighted graph.} which we will denote by $G_{\vD}$.

Vertices of $G_{\vD}$ correspond to the levels of the tiling. The \textit{level} is defined as a horizontal segment, which 
\begin{itemize}[itemsep=-1pt]
\item is equal to the union of several horizontal sides of the tiles; and
\item is maximal, ie, it cannot be extended in either direction while still satisfying the previous condition.
\end{itemize}

Now, we connect the vertices (levels) by the oriented edges corresponding to the tiles, where for every tile $R_i$ we have an edge which ``connects'' two levels containing $R_i$'s top and bottom sides (there could be more than one edge going from vertex $A$ to vertex $B$, and therefore it is quite possible for $G_{\vD}$ to be a multigraph, not a ordinary graph). Every edge is oriented downward---from the higher level to the lower, and every edge $e$ is labeled with the number $\omega_e$ equal to the width of the tile.

\dlpic[16pt]{1.0}{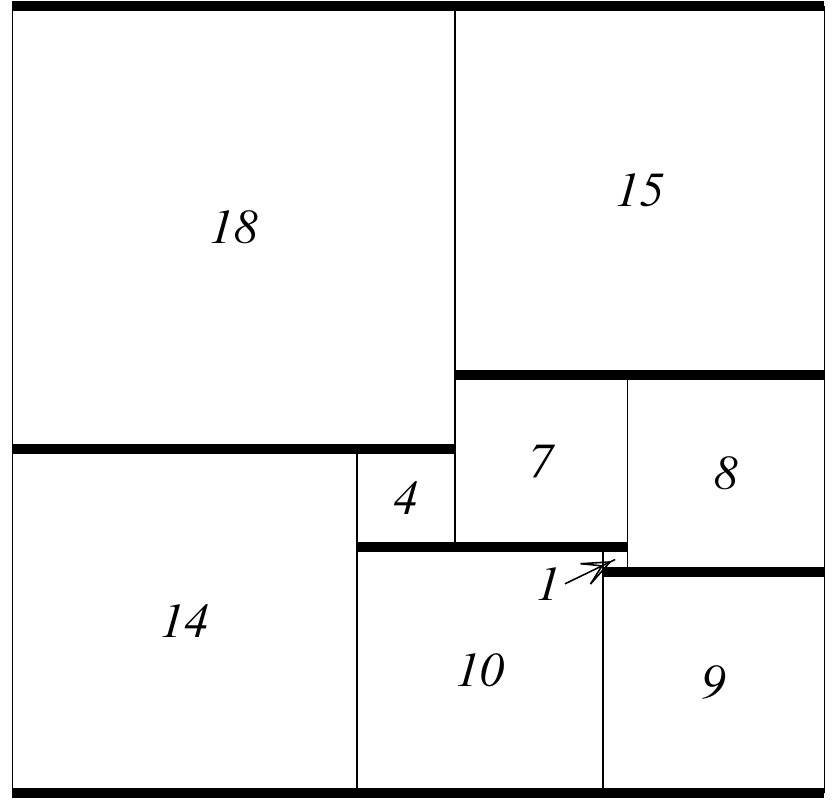}{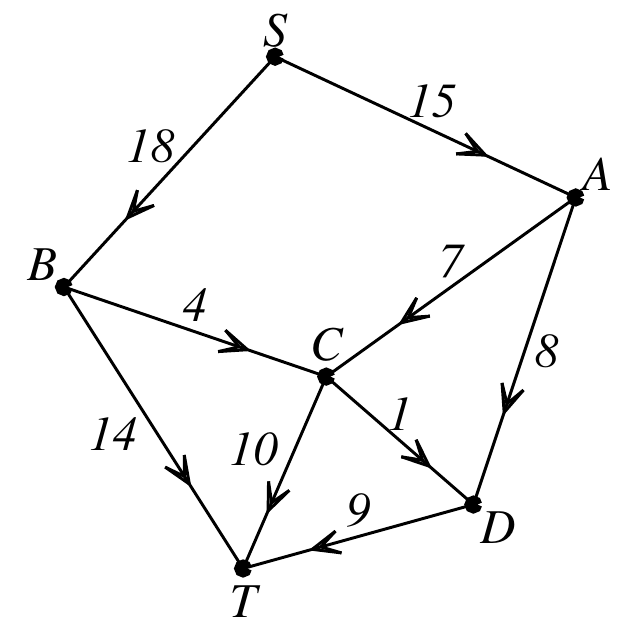}{MoronSqGraph}{A squaring and its network}

For an example of this construction using a well-known squaring of $33\times 32$ rectangle described in \cite{Moron}, see Figure~\ref{fig:MoronSqGraph} on the left. Evidently, there are six levels in this squaring which are shown in \textbf{bold}. The corresponding network is pictured on the right.

Let us denote by $S$ and $T$ the vertices of $G$ which correspond to the top and the bottom sides of $R$, respectively. Also note that graph $G$ is evidently planar and connected. The number of edges in $G$ is the same as the number of tiles $n$ in dissection $\vD$. Thus we are required to prove that $N \lte \tau^n$. Due to Theorem \ref{thm:PlanarGraphSpanTrees} it would suffice to prove the inequality $N \lte t(G) $.

Now, it is easy to see that numbers $\omega_e$ satisfy the system of equations described in Theorem \ref{thm:EdgewiseKirchhoff}. First, it is quite obvious that for any vertex $v$ except for $S$ and $T$ the sum of the weights of all edges leading to $v$ is equal to the sum of the weights of all edges coming out of $v$. Similarly, for any planar face of $G$ the signed sum of the edge weights along its boundary equals zero---we count the weight with $(+)$ sign if its edge is oriented clockwise and with $(-)$ sign otherwise. It follows quite easily that the same is true for any non-self-intersecting cycle in $G$.

Clearly, the throughput of network $G$ is the sum of weights of all edges leading to $T$, which is equal to the length of the bottom side of $R$, ie, to $N$.

Now Theorem \ref{thm:EdgewiseKirchhoff} implies that there exists a unique solution for this system of equations (indeed, the absolute value of its determinant is $t(G) > 0$). Obviously that solution is vector $(\omega_e)$. The right side of that system is vector $(0, 0, \ldots, 0, N)$. From Cramer's rule formulas we have that
$$
\omega_e = \frac {a_e}{t(G)} \,,
$$
where $a_e$ are some integers divisible by $N$, that is, $a_e = N b_e$, $b_e \in \Z$.

Each integer $\omega_e$ is obtained from integer $b_e$ by multiplying it by the same rational fraction $N/t(G)$. Let us assume that $N > t(G)$. Reduce $N/t(G) > 1$ to fraction $N'/m > 1$, where integers $N'$ and $m$ are co-prime. Then numerator $N'$ must be divisible by some prime number $p$, which does not divide denominator $m$. Therefore numbers $\omega_e = b_e (N/t(G)) = b_e (N'/m)$ are all divisible by $p$. This contradicts the assumption that $w_e$ are co-prime as a set, and the lemma is proved.
\end{proof}

Now consider an arbitrary dissection of the given rectangle~$R$ with dimensions $M \times N$.

Assume that one of the non-square tiles, say, $R_k$, has dimensions $a\times b$ with $a > b$. Then perform an ``inverse'' Fibonacci transformation of the pair $(a, b)$, that is, replace $R_k$ with two tiles: $b\times b$ square $K$ and rectangle $R_k'$ with dimensions $(a-b)\times b$.

\clpic{1.0}{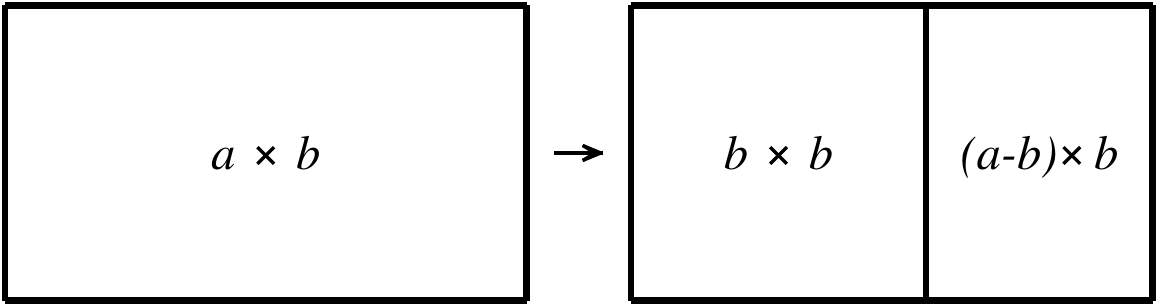}{ifTransform}{``Inverse Fibonacci'' subtiling}

Replacing tiling $\{R_i\}$ with tiling $\{R_i'\} \cup \{K\}$, where for all $i \neq k$ we define $R_i' = R_i$, does not change the right side of inequality \eqref{eqn:sumofdeltas}. At the same time the left side does not increase. Indeed, denote $c = \gcd(a, b)$ and examine the difference $\delta(R_k) - \delta(K) - \delta(R_k')$. We have
$$
\delta(R_k) - \delta(K) - \delta(R_k') = \frac{a}{c} - 1 - \frac{\max(R_k')}{c} = \left( \frac{a}{c} - \frac{\max(R_k')}{c} \right) - 1 \gte 0 \,,
$$
since both ${a}/{c}$ and ${\max(R_k')}/{c}$ are integers with the former greater than the latter, and so their difference must be at least $1$.

Performing a series of such replacements we will eventually transform tiling $\{R_i\}$ into a squaring of $R$. Obviously, all square sizes in this new dissection will still be divisible by $d$. Scaling down by the factor of $d$ we will obtain a prime squaring of rectangle $R'$ with dimensions $(M/d) \times (N/d)$. It should be noted that this operation does not affect neither side of the inequality \eqref{eqn:sumofdeltas}.

The sum of the biases for this squaring does not exceed the sum of the biases for the original tiling, and therefore our theorem follows from Lemma \ref{thm:TauConwayLemma}.
\end{proof}

Theorems \ref{thm:SumOfBiases} and \ref{thm:PlanarGraphSpanTrees} immediately produce two corollaries.

\begin{corollary}
\label{thm:SizeOfPrimeSquaring}
Any prime squaring of integer rectangle $R$ with $\max(R) \gte \tau^n$ consists of at least $n$ tiles.
\end{corollary}

\begin{corollary}
Any squaring of integer rectangle $R$ with $\max(R) \gte \tau^n$, which contains a unit square, consists of at least $n$ tiles.
\end{corollary}

It is rather clear that these bounds are not tight. Incidentally, it was proved in \cite{Trustrum} that the upper bound for the number of squares necessary to tile $N \times N$ square does not exceed $6\log_2(N)$.

\medskip

The last corollary solves items \textit{(a)} and \textit{(b)} of the original problem, proving the lower bound which is slightly better than the one required there. However, the precise answer to item (c) seems to lie somewhere between $1.4 n$ and $3 n$, and finding it would certainly require a very different and, almost surely, much less elementary approach.


\let\cleardoublepage\clearpage

\newcommand\bibauthor[1]{{#1}}
\newcommand\bibyear[1]{\textbf{(#1)}}
\newcommand\bibtitle[1]{\textit{#1}}
\newcommand\bibmisc[1]{{#1}}

\enddocument